\documentclass[12pt]{amsart}
\usepackage{amssymb,amsmath,amsthm}
\newcommand{\leg}[2]{\genfrac{(}{)}{}{}{#1}{#2}}

\newtheorem{theorem}{Theorem}
\newtheorem{lemma}[theorem]{Lemma}
\newtheorem{corollary}[theorem]{Corollary}

\theoremstyle{remark}
\newtheorem*{remark}{Remark}

\newtheorem*{question}{Question}
\numberwithin{theorem}{section} \numberwithin{equation}{section}

\newcommand{\C}{\mathbb{C}}

\newcommand{\Z}{\mathbb{Z}}
\newcommand{\Q}{\mathbb{Q}}

\newcommand{\SL}{{\text {\rm SL}}}

\newcommand{\ve}{\varepsilon}

\newcommand{\Fp}{\mathbb{F}_p}

\newcommand{\ord}{\operatorname{ord}}

\renewcommand{\H}{\mathbb{H}}

\newcommand{\ra}{\rightarrow}

\newcommand{\calH}{\mathcal{H}}

\newcommand{\calO}{\mathcal{O}}

\newcommand{\wtE}{\widetilde{E}}
\newcommand{\wtS}{\widetilde{S}}
\newcommand{\wtx}{\widetilde{x}}

\newcommand{\tr}{\operatorname{tr}}

\newcommand{\End}{\operatorname{End}}
\newcommand{\gen}{\operatorname{gen}}

\newcommand{\Gal}{\operatorname{Gal}}

\newcommand{\9}{^{\phantom9}}

\begin{document}
\title[Reduction of CM elliptic curves and modular function congruences]
{Reduction of CM elliptic curves and modular function
congruences}

\author{Noam Elkies, Ken Ono and Tonghai Yang}

\address{Department of Mathematics, Harvard University,
Cambridge, Massachusetts 02138} \email{elkies@math.harvard.edu}

\address{Department of Mathematics, University of Wisconsin,
Madison, Wisconsin 53706} \email{thyang@math.wisc.edu}
\email{ono@math.wisc.edu}
\thanks{The authors thank the National Science Foundation for its
support. The second author is grateful for the support of the
David and Lucile Packard, H.\ I.\ Romnes and John S. Guggenheim
Fellowships.}
\date{\today}
 \maketitle

\section{Introduction and Statement of Results}

Let $j(z)$ be the modular function for $\SL_2(\Z)$ defined by
$$
j(z):=\sum_{n=1}^{\infty}c(n)q^n=\frac{\left(1+240\sum_{n=1}^{\infty}
\sum_{v\mid n}v^3q^n\right)^3}{q\prod_{n=1}^{\infty}(1-q^n)^{24}}=
q^{-1}+744+196884q+\cdots,
$$
where $q=e^{2\pi i z}$. The coefficients $c(n)$
possess some striking properties. By ``Monstrous Moonshine'',
these integers occur as degrees of a special graded representation
of the Monster group, and they satisfy some classical
Ramanujan-type congruences. In particular, Lehner proved
\cite{lehner} that if $p\leq 7$ is prime and $m$ is a positive
integer, then for every $n\geq 1$ we have
$$
c(p^m n)\equiv \begin{cases}
0\pmod{2^{3m+8}}\ \ \ \  &{\text {\rm if}}\ p=2,\\
0\pmod{3^{2m+3}}\ \ \ \ &{\text {\rm if}}\ p=3,\\
0\pmod{5^{m+1}}\ \ \ \ &{\text {\rm if}}\ p=5,\\
0\pmod{7^m}\ \ \ \ &{\text {\rm if}}\ p=7.
\end{cases}
$$
\noindent Modulo $11$, it also turns out that $c(11n)\equiv
0\pmod{11}$ for every positive integer $n$.

As usual, if $U(p)$ denotes the formal power series
operator
\begin{equation}\label{Up}
\left ( \sum_{n=-\infty}^{\infty} a(n)q^n\right) \ \Big| \ U(p)
:=\sum_{n=
-\infty}^{\infty} a(pn)q^n,
\end{equation}
then these congruences imply that
$$
j(z) \ | \ U(p) \equiv 744\pmod p
$$
for every prime $p\leq 11$. It is natural to ask whether such
congruences hold for any primes $p\geq 13$.

\begin{remark}
Since the Hecke operator $T(p)$ acts like $U(p)$ on spaces of
holomorphic integer weight modular forms modulo $p$, this problem
is somewhat analogous to that of determining whether there are infinitely
many non-ordinary primes for the generic integer weight newform
without complex multiplication.  Apart from those newforms
associated to modular elliptic curves, for which
the existence of infinitely many non-ordinary primes
was shown in~\cite{Elkies}, little is known.
\end{remark}

Serre showed \cite{serre} that the answer to this question
for $j(z)$ is
negative, an observation which has been
generalized by the second author and Ahlgren \cite{compositio}.
In particular,  if $F(x)\in \Z[x]$ is a polynomial
of degree $m\geq 1$ and $p > 12m+1$ is a prime which does not
divide the leading coefficient of $F(x)$, then (see Corollary 5 of
\cite{compositio})
\begin{displaymath}
  F(j(z)) \ | \ U(p) \not \equiv a(0)\pmod p,
\end{displaymath}
where $a(0)$ is the constant term in the Fourier expansion of
$F(j(z))$.

\begin{remark}
In the cuspidal case where
$(m,p)=(1,13)$, it is interesting to note that
(for example, see Section (6.16)
or \cite{serre})
$$
(j(z)-744)\ | \ U(13)\equiv -\Delta(z)\pmod{13},
$$
where $\Delta(z)=q\prod_{n=1}^{\infty}(1-q^n)^{24}$
is the usual Delta-function.
\end{remark}

Here we investigate the more general question concerning the
existence of congruences of the form
\begin{equation}\label{congtype}
F(j(z))\ | \ U(p) \equiv G_p(j(z))\pmod p,
\end{equation}
where $G_p(x)\in \Z[x]$.  The result quoted above implies
that congruences of the form (\ref{congtype}) do not hold for any
primes $p>12m+1$. This follows from the simple fact that the only
polynomials in $j(z)$ whose Fourier expansions do not contain
negative powers of $q$ are constant.

In Section~\ref{sect:2} we give
a general criterion (Theorem~\ref{Uptheorem})
that proves such congruences, and we apply it to Hilbert class
polynomials.  For a discriminant $-D<0$, let
$\calH_D(x)\in \Z[x]$ be the associated Hilbert class polynomial.
More precisely, $\calH_D(x)$ is the polynomial of degree $h(-D)$
whose roots are the singular moduli of discriminant~$-D$,
where $h(-D)$ is the class number of the ring of integers
$\calO_D$ of $\Q(\sqrt{-D})$.
By the theory of complex multiplication, these singular moduli
are the $j$-invariants of those elliptic curves that have complex
multiplication by $\mathcal{O}_D$, the ring of integers of
$\Q(\sqrt{-D})$.

Although there are no congruences of the form (\ref{congtype}) for
$F(x)=\calH_D(x)$ involving primes $p> 12h(-D)+1$, we show that
such congruences are quite common for smaller primes, as Lehner
demonstrated for $F(x)=\calH_{3}(x)=x$.

\begin{theorem}\label{congtheorem}
Suppose that $-D<0$ is a fundamental discriminant, and that
integers $c_D(n)$ are defined by
$$
\calH_D(j(z))=\sum_{n=-h(-D)}^{\infty} c_D(n)q^n.
$$
\begin{enumerate}
\item If $p\leq 11$ is prime and $h(-D)<p$, then
$$
\calH_D(j(z))\ | \ U(p)\equiv c_D(0)\pmod p.
$$

\item For every prime~$p$ there is a non-positive integer $N_p$
such that
$$
\calH_D(j(z))\ | \ U(p)\equiv G_p(j(z))\pmod p,
$$
for some $G_p(x)\in \Z[x]$,
provided that $-D<N_p$ and $\leg{-D}{p}\neq 1$.
This polynomial $G_p(x)$ has degree $\leq h(-D)/p$.
In particular, $G_p(x)$ is constant if $h(-D)<p$,
or more generally if $c_D(-pn)\equiv 0\pmod p$ for every integer $n>0$.
\end{enumerate}
\end{theorem}

Proving Theorem~\ref{congtheorem}~(2) depends heavily on the
interplay between singular moduli and supersingular
$j$-invariants.  For a prime $p\geq 5$, define the supersingular loci
$S_p(x)$ and $\wtS_p(x)$ in $\Fp[x]$ by
\begin{equation}\label{sslocus}
\begin{split}
S_p(x)&:= \prod_{E/\overline{\mathbb{F}}_p\ {\text {\rm
supersingular}}} (x-j(E)),\\
\wtS_p(x)&:= \prod_{\substack{E/\overline{\mathbb{F}}_p \
\ {\text {\rm supersingular}}\\ j(E) \not \in \{0, 1728\}}}
(x-j(E)).
\end{split}
\end{equation}
These products are over isomorphism classes of supersingular
elliptic curves. It is a classical fact (for example, see
\cite{Si}) that the degree of $\wtS_p(x)$ is $\lfloor
p/12\rfloor$.

The criterion for proving
congruences of the form (\ref{congtype}) is stated in terms of the
divisibility of $F(x)$ by $\wtS_p(x)^2$ in $\Fp[x]$. For
Hilbert class polynomials, this criterion is quite natural  since
a classical theorem of Deuring asserts that the reduction of every
discriminant $-D$ singular modulus modulo $p$
is a supersingular $j$-invariant when
$p$ does not split in $\Q(\sqrt{-D})$.

Therefore, to prove Theorem~\ref{congtheorem} we are forced to
consider the surjectivity of Deuring's reduction map of singular
moduli onto supersingular $j$-invariants, a question which is
already of significant interest.
Using classical facts about elliptic curves with CM and
certain quaternion algebras, we reinterpret this problem in terms of the
vanishing of Fourier coefficients of
specific weight $3/2$ theta functions constructed by Gross.
Then, using deep results of Duke and Iwaniec which bound
coefficients of half-integral weight cusp forms, we obtain the
following theorem.

\begin{theorem}\label{surjectivity}
If $p$ is an odd prime and $t\geq 1$, then there is a
non-positive integer $N_p(t)$ such that $S_p(x)^t \mid \calH_D(x)$
in $\Fp[x]$ for every fundamental discriminant $-D<N_p(t)$ for
which $\leg{-D}{p}\neq 1$.
\end{theorem}

\noindent {\it Three Remarks.} 1. After this paper was submitted
for publication, Bill Duke and the referee informed us of earlier
work of P. Michel \cite{michel}. In this recent paper, Michel
obtains equidistribution results which imply Theorem
\ref{surjectivity} for discriminants for which $p$ is inert. His
proof, which is different from ours, is based on subconvexity
bounds for $L$-functions. Our proof, which also includes the
ramified cases, is based on non-trivial estimates of Fourier
coefficients of half-integral weight cusp forms. Both proofs are
somewhat related via Waldspurger's formulas connecting values of
$L$-functions to Fourier coefficients.

2. Theorem \ref{surjectivity} is ineffective due to the
ineffectivity of Siegel's lower bound for class numbers. It is
possible to obtain effective results by employing various Riemann
hypotheses (for example, see work \cite{OS} by the second author
and Soundararajan concerning Ramanujan's ternary quadratic form),
or by assuming the non-existence of Landau-Siegel zeros.

3. Theorem~\ref{surjectivity} is closely related to the work of
Gross and Zagier \cite{GZ} which provides the prime factorization
of norms of differences of singular moduli in many cases.

\medskip

In Section~\ref{sect:2} we use a result of Koike (arising in his study
of ``$p$-adic rigidity of $j(z)$'') to prove Theorem~\ref{Uptheorem}.
In Section~\ref{sect:3} we recall preliminary facts regarding
endomorphism rings of elliptic curves with complex multiplication
and quaternion algebras, and  we prove Theorem~\ref{surjectivity}
using facts about weight $3/2$ Eisenstein series combined with the
Duke-Iwaniec bounds for coefficients of weight $3/2$ cusp forms.
Then we combine these results with Theorem~\ref{Uptheorem} to
prove Theorem~\ref{congtheorem}. In Section~\ref{sect:4} we
conclude with some remarks on numerical calculations related to
Theorems~\ref{congtheorem} and \ref{surjectivity}.

\section{A congruence criterion and supersingular $j$-invariants}\label{sect:2}

Here we give a simple criterion which implies congruences of the
form (\ref{congtype}). This criterion is a simple generalization
of Theorem~2 of \cite{compositio}.
The following result of Koike (see Proposition~1 of \cite{Koike}),
which is a special case of work of Dwork and Deligne \cite{dwork},
describes the Fourier expansion of $j(pz)\pmod{p^2}$
in terms of $j(z)$ and the collection of supersingular $j$-invariants.
Since clearly $j(pz) \equiv j(z)^p \bmod p$,
we can describe the Fourier expansion of $j(pz)\pmod{p^2}$
via the reduction modulo~$p$ of the expansion of $(j(pz) - j(z)^p) / p$.

\begin{theorem}\label{koike}
For each prime $p$ there is a rational function
$\delta_p(x) = N_p(x) / \wtS_p(x)$, with $N_p \in \Fp[x]$,
such that
\begin{equation}
j(pz) \equiv j(z)^p + p \delta_p(j(z)) \pmod{p^2}.
\end{equation}
\end{theorem}

\begin{corollary}\label{koike_cor}
For all $F \in \Z[x]$ we have
\begin{equation}
F(j(pz)) \equiv F(j(z)^p) + p F'(j(z)^p) \delta_p(j(z)) \pmod{p^2},
\end{equation}
where $F'$ denotes the derivative of~$F$.
\end{corollary}

\begin{proof}
By linearity it suffices to prove this for $F=x^k$ ($k=0,1,2,\ldots$).
The case $k=0$ is trivial; $k=1$ is Theorem~\ref{koike};
and for $k>1$ we may either raise the congruence in Theorem~\ref{koike}
to the power~$k$\/ and reduce mod~$p^2$, or argue by induction from
the case $k-1$ and the same congruence, obtaining
\begin{equation}
j(pz)^k \equiv j(z)^{pk} + p j(z)^{(k-1)p} \delta_p(j(z)) \pmod{p^2}
\end{equation}
as claimed.
\end{proof}

\noindent We can now deduce our congruence criterion
using little more than the theory of Hecke operators.

\begin{theorem}\label{Uptheorem}
Let $F(x)\in \Z[x]$ be a monic polynomial of degree~$m$, and let
$$
F(j(z))=\sum_{n=-m}^{\infty}a(n)q^n.
$$
\begin{enumerate}
\item If $p$ is prime and $\wtS_p(x)^2$
divides $F(x)$ in $\Fp[x]$, then
$$
F(j(z)) \ |\ U(p)\equiv G_p(j(z))\pmod p,
$$
for some $G_p(x)\in \Z[x]$ with degree $\leq m/p$.
\item If $p\leq 11$ is prime and $m<p$, then
$$
 F(j(z))\ | \ U(p)\equiv a(0)\pmod p,
$$
where $a(0)$ is the constant term in the Fourier expansion of
$F(j(z))$.
\end{enumerate}
\end{theorem}

\begin{proof}
For (1), let
$$
F(j(z))=\sum_{n= -\infty}^{\infty} a(n)q^n.
$$
Denote by $T_0(p)$ the operator $pT(p)$, that is,
$p$ times the usual weight zero $p$th Hecke operator.
Then we have
\begin{equation}\label{sstag5.3}
  p F(j(z)) \ |\ U(p)
= F(j(z)) \ |\ T_0(p) -  F(j(pz))
= p\sum_{n= -\infty}^{\infty}a(pn)q^n.
\end{equation}
The modular function $F(j(z))\ |\ T_0(p)$ is in $\Z[j(z)]$ since
it has integer Fourier coefficients and is holomorphic on $\H$.
Since $\wtS_p(x)^2 | F(x)$ we have
$\wtS_p(x) | F'(x)$,
whence also $\wtS_p(x) | F'(x^p)$.
By Corollary~\ref{koike_cor} it follows that
$F(j(pz))\pmod{p^2}$ is congruent modulo~$p^2$
to an integer polynomial in $j(z)$,
namely $F(j(z))^p + p F'(j(z)^p) \delta_p(j(z))$.
Thus (\ref{sstag5.3}) yields the desired congruence modulo~$p$
between $F(j(z))\ |\ U(p)$ and a polynomial in $j(z)$.
Moreover, this polynomial must have degree $\leq m/p$
because $F(j(z)) \ |\ U(p)$ has valuation $\geq -m/p$ at the cusp.

For (2), we observe that the condition $\wtS_p(x)^2 | F(x)$
of~(1) is vacuous for $p\leq 11$, because $\wtS_p(x)=1$ for those~$p$.
Thus $F(j(z)) \ | \ U(p)$ is always congruent mod~$p$ to a polynomial
in~$j$ of degree at most $m/p$.  In particular, if $m<p$ then this
polynomial reduces to a constant, which must equal $a(0)$
by the definition of $U(p)$.

\end{proof}

\section{Gross' theta functions
and the proofs of Theorems~\ref{congtheorem} and
\ref{surjectivity}}\label{sect:3}

Throughout, $p$ shall denote a prime. We begin by recalling
certain facts about elliptic curves with complex multiplication
(for example, see Chapter II of \cite{silverman}).
Let $B$ be the unique quaternion algebra over~$\Q$ ramified
exactly at $p$ and~$\infty$.  For $x \in B$,
let $Q(x) = x \bar x = -x^2$, the reduced norm of~$x$;
the map $Q: B \ra \Q$ is a quadratic form on~$B$,
which is positive-definite because $B$\/ is ramified at~$\infty$.
Fix a maximal order $R \subset B$.
Then $Q$\/ takes integer values on~$R$,
and since $B$\/ is ramified at~$p$, the subset
$$
\pi := \{ x \in R  \mid  p|Q(x) \}
$$
of~$R$\/ is a two-sided ideal with $R/\pi$
a finite field of $p^2$ elements and $\pi^2 = pR$.

Let $K=\Q(\sqrt{-D})$ be an imaginary quadratic field
with ring of integers $\mathcal O_D$.
More generally, for a positive integer $m$ congruent to~$0$ or~$3$
mod~$4$, let $$\mathcal O_m = \Z + \frac12(m+\sqrt{-m}) \Z,$$
the order of discriminant $-m$ in $\Q(\sqrt{-m})$.
An {\em optimal embedding} of $\calO_m$ into~$R$ is an embedding in
$i: \Q(\sqrt{-m}) \hookrightarrow B$ for which $i^{-1}(R) = \calO_m$.
Such an embedding is determined by the image of $\sqrt{-m}$ in
$$
V =\{ x \in B \mid \tr x=0 \},
$$
a 3-dimensional subspace of~$B$.  This image must be an element
of norm~$m$\/ in the lattice
$$
L := V \cap (\mathbb Z + 2 R)
$$
in~$V$, and conversely every $v\in L$ of norm~$m$
comes from an embedding $\Q(\sqrt{-m}) \hookrightarrow B$,
which is optimal if and only if $v$ is a primitive vector of~$L$
(that is, $v \notin fL$ for any $f>1$).

Two optimal embeddings $i_1$, $i_2$ are {\em equivalent}\/
if they are conjugate to each other by a unit in~$R$\/;
In other words, if there is $u \in R^{\times}$ such that
$
i_1(x) = u i_2(x) u^{-1}
$
for all $x \in \calO_m$. Let $h(\calO_m, R)$ be the number of
equivalence classes of optimal embeddings of $\calO_m$ into $R$.
Using the connection between embeddings and lattice vectors,
Gross proved \cite{Gross} that these numbers generate
the theta series of~$L$ as follows:

\begin{lemma} \label{lemma3.1}{\text {\rm {(\cite{Gross}, Proposition 12.9)}}}
The theta function
$$
\theta_L(z) := \sum_{x \in L} e(Q(x))
$$
is given by
$$
\theta_L(z) = 1 + \sum_{m \ge 1} a_R\9(m)q^m,
$$
where
$$
a_R\9(m) = w_R\9 \sum_{m = f^2 D} \frac{h(\calO_D, R)}{u(d)},
$$
in which $u(D)=\frac12 \# \calO_d^{\times}$ and $w_R\9 = \# R^{\times}$.
Moreover, $\theta_L$ is a holomorphic modular form
of weight $3/2$ and level $4p$.
\end{lemma}
\begin{remark}
This $\theta_L(z)$ is a modular form of half-integral weight
in the sense of Shimura~\cite{shimura}.
Moreover, it lies in Kohnen's plus-space \cite{K}.
\end{remark}

Now recall that every maximal order of~$B$\/ is isomorphic with
the endomorphism ring of some supersingular elliptic curve $E_0$
over $\overline{k}=\overline{\mathbb{F}}_p$, say $R=\End(E_0)$,
with the ideal~$\pi$ comprising the inseparable endomorphisms of~$E_0$.
Let $\mathbb C_p$ be the completion of an algebraic closure of~$\Q_p$.
The residue field of the unramified quadratic extension of~$\Q_p$
in~$\mathbb \C_p$ is a finite field of~$p^2$ elements; call it~$k$.
There is then a canonical map $R \ra k$, $x \mapsto \wtx$,
defined as follows: any $x \in R = \End(E_0)$ induces
multiplication by $\wtx$ on the invariant differentials of~$E_0$.
The kernel $\{x \mid \wtx = 0\}$ is our ideal~$\pi$,
so this map $x \mapsto \wtx$ identifies~$k$ with $R/\pi$.

Suppose that $(E, i)$ is a CM elliptic curve over $\mathbb C_p$
with complex multiplication by~$\calO_D$.  We may then choose a map
$$
i: \calO_D \stackrel{\sim}{\ra} \End(E)
$$
which is {\em normalized}\/ in the sense that any $a\in\calO_D$
acts on the invariant differentials of~$E$\/ by multiplication by~$a$.
(There are two choices of~$i$, related by conjugation in $\Gal(K/\Q)$,
and one of them is normalized.)
If $p$ is inert or ramified in~$K$, then a classical result of Deuring
states that $\wtE := E\mod \mathfrak p$ is a supersingular
elliptic curve over $\overline{k}$, and if $\wtE \cong E_0$,
then we obtain an optimal embedding
$$
f: \calO_D \cong \End(E) \ra \End(\wtE)
\cong \End(E_0)=R,
$$
and moreover the embedding is {\em normalized}\/:
if $x=f(a)$ then $\wtx$ is the residue of~$a$ in~$k$.
Any embedding equivalent to a normalized one is again normalized,
because conjugation by a unit in $R^{\times}$
commutes with our map $x \mapsto \wtx$.
Since any two isomorphisms $\wtE \cong E_0$ differ by
multiplication by a unit in $R^{\times}$, one sees that
the equivalence class of~$f$\/ is uniquely determined by $(E,i)$.
Conversely, given a normalized optimal embedding $f: \calO_D \ra R$,
Deuring's lifting theorem (see \cite[Proposition 2.7]{GZ})
asserts that there is a CM elliptic curve $(E, i)$,
unique up to isomorphism, such that $\wtE \cong E_0$, and
its associated optimal embedding is normalized and equivalent to~$f$.

Note now that embeddings of~$\calO_K$ into~$R$\/
come in conjugate pairs $\{i,\bar\imath\}$, where
$\bar\imath(a) = i(\bar a) = \overline{i(a)}$.
If $p$ is inert in~$K$,\/ then in each pair $\{i,\bar\imath\}$
exactly one of $i,\bar\imath$ is normalized, whereas
if $p$ is ramified in~$K$\/ then every embedding is normalized.
We have thus proved the following lemma.

\begin{lemma}\label{lemma3.2} Let $J_D$ be the set of $j$-invariants of
CM elliptic curves with endomorphism ring $\calO_D$.
If we set
$$
J_D(E_0) = \{ j \in J_D \mid j \bmod \mathfrak p = j(E_0)\},
$$
then
$$
\# J_D(E_0) = \ve h(\calO_D, R),
$$
where $R=\End(E_0)$ and $\ve=1/2$ or~$1$
according as $p$ is inert or ramified in~$K$.
\end{lemma}

\begin{proof}[Proof of Theorem~\ref{surjectivity}]
For every prime $p$, there are finitely many supersingular
elliptic curves $E_0$ over $\overline{\mathbb{F}}_p$. So to prove
the theorem it suffices to show, for each supersingular curve $E_0$
with $R=\End(E_0)$ and each positive integer $t$, that
there is a non-positive integer $N_p(t)$ such that  every
fundamental discriminant $-D < N_p(t)$ with $\leg{-D}{p} \ne 1$
has the property that
$\ord_{x=j(E_0)}(\calH_D(x) \bmod \mathfrak p) \ge t$. By
Lemmas \ref{lemma3.1} and \ref{lemma3.2}, this is equivalent to
 \begin{equation}
 \label{inequality}
 a_R\9(D) \ge \frac{w_R\9 t}{\ve u(d)}.
\end{equation}

This turns out to be a simple consequence of well-known deep
results of Siegel \cite{siegel}, Duke \cite{Duke}, and
Iwaniec \cite{iwaniec}. Indeed, one has by \cite{siegel}
\begin{equation}
\theta_L(z) =\frac{12}{p-1}E_{\gen(L)}(z) + C_L(z),
\end{equation}
where $E_{\gen(L)}(z)$ is the Eisenstein series associated
to the genus of the lattice $L$ and $C_L(z)$ is a cusp
form of weight $3/2$ and level $4p$.  Although Siegel's result
is not stated for forms of half-integral weight forms,
the proof follows {\it mutatis mutandis}, with the constant
$12/(p-1)$ coming from Gross' explicit calculation
of these Eisenstein series.
More precisely, he shows in \cite[(12.11)]{Gross}
(this is $2G$ in his notation, see also \cite[\S8]{KRY}) that
\begin{equation}
E_{\gen(L)}(z)
 =\frac{p-1}{12} + 2 \sum_{m>0} H_p(m) q^m,
\end{equation}
where $H_p(m)$ is a slight modification of Kronecker-Hurwitz
class number $H(m)$ defined in \cite[(1.8)]{Gross}.
In particular, when $-D$ is a fundamental discriminant, one has
$$
H_p(D) =\frac{1}2 \left( 1 -\leg{-D}{p}\right ) \frac{h(-D)}{u(D)},
$$
where $h(-D)$ is the ideal class number of $K$.

By Siegel's theorem \cite{siegel}, one has
$$
H_p(D) \gg  D^{\frac{1}2- \epsilon}
$$
for every $\epsilon >0$ when $ (\frac{-D}p)\ne 1$. On the other
hand, a theorem of Duke  \cite{Duke}, which extended earlier work
of Iwaniec \cite{iwaniec}, implies that the coefficients of the
cusp form
$$
C_L(z) =\sum_{m \ge 1} c_L(m) q^m
$$
satisfy
$$
|c_L(D)| \ll D^{\frac{3}{7}+\epsilon}.
$$
Since $3/7 < 1/2$, we get
$$
a_R\9(D) = \frac{24}{p-1}H_p(D) + c_L(D) \gg D^{\frac{1}2 -\epsilon}.
$$
This proves (\ref{inequality}), and consequently completes the
proof of the theorem.
\end{proof}

\begin{proof}[Proof of Theorem~\ref{congtheorem}]
Theorem~\ref{congtheorem}~(1) follows immediately from
Theorem~\ref{Uptheorem}~(1).  Theorem~\ref{congtheorem}~(2)
follows from Theorem~\ref{Uptheorem}~(2) and
Theorem~\ref{surjectivity} by letting $N_p=N_{p}(2)$.
\end{proof}

\section{Concluding Remarks}\label{sect:4}

Numerical computations reveal many nearly uniform sets of examples
of congruences of the form $(\ref{congtype})$. Here we comment on
those cases where
$$
\calH_D(j(z))\ | \ U(p)\equiv c_D(0)\pmod p.
$$
In view of Theorem~\ref{Uptheorem}, it is natural to investigate
those fundamental discriminants $-D<0$ for which
\begin{equation}\label{range}
p/6 < h(-D) < p.
\end{equation}
The lower bound of this inequality is dictated by the fact that the
degree of $\wtS_p(x)$ is $\lfloor p/12\rfloor$, and the
upper bound is chosen so that the $U(p)$ operator does not produce
a Fourier expansion with negative powers of $q$.

Computations reveal
that if $-239 < -D < 0$ and $p$ is an odd prime satisfying
(\ref{range}) for which $\leg{-D}{p}\neq 1$, then $S_p(x)^2$
divides $\calH_D(x)$ in $\Fp[x]$, which, by
Theorem~\ref{Uptheorem}, in turn implies that
$$
   \calH_D(j(z)) \ | \ U(p)\equiv c_D(0)\pmod p.
$$
This uniformity suggests that this phenomenon might hold in
generality. However, this is not true; when $-D=-239$, we have
$$
\calH_{-239}(j(z))\ | \ U(79) \equiv 44+2q+62q^2+\cdots\pmod{79},
$$
although $79$ is inert in $\Q(\sqrt{-239})$ and $h(-239)=15$.  In
this case $S_{79}(x)$ divides $\calH_{239}(x)$ in
$\mathbb{F}_{79}[x]$, but the supersingular $j$-invariant $j=-15$ is
a root of multiplicity only $1$.
This raises the following natural question.

\begin{question}
If $p$ is an odd prime, then define $\Omega_p$ by
$$
\Omega_p:=\left \{ -D\ {\text {\rm fundamental}} \; \big| \;
p/6 < h(-D) < p\ {\text {\rm and}}\ \leg{-D}{p}\neq 1\right \}.
$$
In general, what ``proportion'' of $-D\in \Omega_p$ have the property
that
$$
\calH_D(j(z))\ | \ U(p) \equiv c_D(0)\pmod p?
$$
\end{question}


\begin{thebibliography}{999999}

\bibitem[AO]{compositio} S.~Ahlgren and K.~Ono,
\emph{Arithmetic of singular moduli and class equations},
Compositio Mathematica  \textbf{141} (2005), pages 293--312.

\bibitem[Du]{Duke} W.~Duke,
\emph{Hyperbolic distribution problems and half-integral weight Maass forms},
Invent.\ Math.\ \textbf{92} (1988), pages 73--90.

\bibitem[Dw]{dwork} B.~Dwork, \emph{$p$-adic cycles},
Inst.\ Hautes \'Etudes Sci. Publ. Math. \textbf{37} (1969),
pages 27--115.

\bibitem[El]{Elkies} Elkies, N.D.:
\emph{The existence of infinitely many supersingular primes
for every elliptic curve over~$\Q$},
Invent.\ Math.\ \textbf{89} (1987), pages 561--567.

\bibitem[Gr]{Gross} B.~Gross,
\emph{Heights and the special values of L-series},
Number Theory (Montreal, Quebec, 1985) CMS Conference Proc.\ \textbf{7},
Amer.\ Math.\ Soc.\ (1987), pages 115--187.

\bibitem[GZ]{GZ} B.~Gross and D.~Zagier,
\emph{On singular moduli},
J.\ reine angew.\ Math.\ \textbf{355} (1985), pages 191--220.

\bibitem[Iw]{iwaniec} H.~Iwaniec,
\emph{Fourier coefficients of modular forms of half-integral weight},
Invent.\ Math.\ \textbf{87} (1987), pages 385--401.

\bibitem[Koh]{K} W.~Kohnen,
\emph{Newforms of half-integral weight},
J.\ reine angew.\ Math.\ \textbf{333} (1982), pages 32--72.

\bibitem[Koi]{Koike} M.~Koike,
\emph{Congruences between modular forms and functions
and applications to the conjecture of Atkin},
J.\ Fac.\ Sc.\ Univ.\ Tokyo, Sect.\ IA Math.\ \textbf{20} (1973),
pages 129--169.

\bibitem[KRY]{KRY} S.~Kudla, M.~Rapoport, and T.H.~Yang,
\emph{The derivative of Eisenstein series and the Faltings's heights},
Compos.\ Math.\ \textbf{ 140} (2004), pages 887--951.

\bibitem[Le]{lehner} J.~Lehner,
\emph{Further congruence properties of the Fourier coefficients
of the modular invariant $j(\tau)$},
Amer.\ J.\ Math.\ \textbf{71} (1949), pages 373--386.

\bibitem[Mi]{michel} P. Michel,
\emph{The subconvexity problem for Rankin-Selberg $L$-functions
and equidistribution of Heegner points}, Annals of Math.
\textbf{160} (2004), pages 185--236.

\bibitem[OS]{OS} K.~Ono and K.~Soundararajan,
\emph{Ramanujan's ternary quadratic form,}
Invent.\ Math.\ \textbf{130} (1997), pages 415--454.

\bibitem[Se]{serre} J.-P.~Serre,
\emph{Divisibilit\'e de certaines fonctions arithm\'etiques},
L'Enseign.\ Math.\ \textbf{22} (1976), pages 227--260.

\bibitem[Sh]{shimura} G.~Shimura,
\emph{On modular forms of half-integral weight},
Ann.\ of Math.\ \textbf{97} (1973), pages 440--481.

\bibitem[Si]{siegel}  C.L.~Siegel,
\emph{\"Uber die analytische Theorie der quadratischen Formen I, II, III},
Ann.\ of Math.\ \textbf{36} (1935), pages 527--606;
\textbf{37} (1936), pages 230--263;
\textbf{38} (1937), pages 212--291.

\bibitem[Si1]{Si} J.\ Silverman,
\emph{The Arithmetic of Elliptic curves},
Springer-$\!$Verlag, New York, 1986.

\bibitem[Si2]{silverman} J.\ Silverman,
\emph{Advanced Topics in the Arithmetic of Elliptic Curves},
Springer-$\!$Verlag, New York, 1994.


\end{thebibliography}
\end{document}